# A MULTIVARIATE EMPIRICAL BAYES STATISTIC FOR REPLICATED MICROARRAY TIME COURSE DATA

By Yu Chuan Tai[1] and Terence P. Speed[2]

*University of California, Berkeley and Walter and Eliza Hall Institute of Medical Research, Australia*

In this paper we derive one- and two-sample multivariate empirical Bayes statistics (the *MB*-statistics) to rank genes in order of interest from longitudinal replicated developmental microarray time course experiments. We first use conjugate priors to develop our one-sample multivariate empirical Bayes framework for the null hypothesis that the expected temporal profile stays at 0. This leads to our one-sample *MB*-statistic and a one-sample $\widetilde{T}^2$-statistic, a variant of the one-sample Hotelling $T^2$-statistic. Both the *MB*-statistic and $\widetilde{T}^2$-statistic can be used to rank genes in the order of evidence of nonzero mean, incorporating the correlation structure across time points, moderation and replication. We also derive the corresponding *MB*-statistics and $\widetilde{T}^2$-statistics for the one-sample problem where the null hypothesis states that the expected temporal profile is constant, and for the two-sample problem where the null hypothesis is that two expected temporal profiles are the same.

**1. Introduction.** Microarray time course experiments differ from other microarray experiments in that gene expression values at different time points can be correlated. This may happen when the design is longitudinal, that is, where the mRNA samples at successive time points are taken from the same *units*. Such longitudinal experiments make it possible to monitor and study the temporal *changes* within units of biological processes of interest for thousands of genes simultaneously. Two major categories of time course experiments are those involving periodic and developmental phenomena. Periodic time courses typically concern natural biological processes such

Received August 2004; revised August 2005.

[1]Supported in part by NIH Grant HG00047.

[2]Supported in part by NIH Grant R01 Lmo7609-01.

as cell cycles or circadian rhythms, where the temporal profiles follow regular patterns [7, 8, 25, 26]. On the other hand, developmental time course experiments measure gene expression levels at a series of times in a developmental process, or after applying a treatment such as a drug to the organism, tissue or cells [9, 32, 34]. In this case, we typically have few prior expectations concerning the temporal patterns of gene expression. The gene ranking methods we develop in this paper are mainly for longitudinal replicated developmental time course data.

A typical microarray time course dataset consists of expression measurements of $G$ genes across $k$ time points, under one or more biological conditions (e.g., wildtype versus mutant). The number of genes $G$ (1,000–40,000) is very much larger than the number of time points $k$, which can be 5–10 for shorter and 11–20 for longer time courses. Many such experiments are unreplicated due to cost or other limitations, and when replicates are done, the number $n$ is typically quite small, say, 2–5. We refer the reader to [29] for a fuller review of microarray time course experiments.

One of the statistical challenges here is to identify genes of interest. In what we call the one-sample problem, these are genes whose patterns of expression change over time, perhaps in some specific way. In the two-sample problem we seek genes whose temporal patterns differ across two biological conditions. Such genes are of interest to biologists because they are often involved in the biological processes motivating the experiment. The challenge arises from the fact that there are very few time points, and very few replicates per gene. The series are usually so short that we cannot consider using standard time series methods as described in [10], such as Fourier analysis, ARMA models or wavelets. The methods proposed in this paper are for the one- and two-sample problems with longitudinal replicated microarray time course experiments of the developmental kind.

The gene ranking problem for such microarray experiments is relatively new. Few methods have been proposed to deal specifically with these problems. The most widely used method for identifying temporally changing genes in replicated microarray experiments is to carry out multiple pairwise comparisons across times, using statistics developed for comparing two independent samples, [2, 6, 13, 19, 20, 23, 24, 33]. These methods are not entirely appropriate as they do not incorporate the fact that longitudinal microarray time course samples are correlated. A simple and intuitive approach to our problem is to use classical or mixed ANOVA models; see Chapter 6 of [11] for a discussion of the latter for analyzing longitudinal data and [22] for a modified approach based on the former for use in the microarray context. However, a number of questions are not adequately addressed by the classical ANOVA methods, or the variants of [22]. As with the pairwise comparisons, the $F$-statistic assumes that gene expression measurements at different times are independent. The classical ANOVA models also assume normality of the gene expression measurements, which may not be a great concern when these are on the log scale. More importantly, standard $F$-statistics in this context will generally lead to more false positives and false negatives than is desirable, due to poorly estimated error variances in the denominator. This issue can be dealt with using the idea of moderation; see, for example, [2, 14, 20, 24, 33]. Moderation in our longitudinal context means the smoothing of gene-specific sample variance–covariance matrices toward a common matrix. When we do this, fewer genes which are not differentially expressed over time but have very small replicate variances are falsely identified as being differentially expressed, and fewer genes which are differentially expressed over time but have large replicate variances (e.g., due to outliers) are missed by the $F$-statistic.

A moderated gene-specific score based on the Wald statistic for the longitudinal one-sample problem was proposed in [15]. However, their method is only applicable to the situation when the number of replicates is greater than the number of time points. In [3] the expression profiles for each gene and each of two biological conditions were represented by continuous curves fitted using B-splines. A global difference between the two continuous curves and an ad hoc likelihood-based $p$-value was calculated for each gene. B-splines were also used in [17] to identify genes with different temporal profiles in the two-sample case. Recently, B-splines were again adopted by [27] to model the population mean, constructing the $F$-statistics for both longitudinal and cross-sectional data with one or more biological conditions. A major feature of this paper was a careful treatment of the multiple testing issue in this context. A novel HMM approach which incorporates the dependency in gene expression measurements across times was proposed in [36] for data with two or more biological conditions. This is one example of using HMM to identify differentially expressed genes across at least two biological conditions in this context.

The multivariate empirical Bayes model proposed in this paper was motivated by the analogous univariate model proposed in [20] for identifying differentially expressed genes in two-color comparative microarray experiments, and the more recent extensions by Smyth [24]. The $B$-statistic in [20] and [24], and the univariate moderated $t$-statistic $\tilde{t}_g$ in [24], consider just one parameter or contrast at a time in the null hypotheses. They are not for null hypotheses with two or more parameters or contrasts of interest simultaneously. However, a partly-moderated $F$-statistic was introduced in [24], which moderates the error variance in the denominator of the ordinary $F$. This partly-moderated $F$-statistic is useful for the simultaneous comparison of multiple uncorrelated contrasts, but as mentioned above, it is not appropriate for longitudinal experiments. Both the $MB$-statistics and the $\tilde{T}^2$-statistics derived in this paper provide a degree of moderation, while retaining the temporal correlation structure. They can be used with both single- and two-channel microarray experiments.

This paper is organized as follows. After briefly introducing our notation, we formally state the null and alternative hypotheses for the gene ranking problem in Section 2. In Section 3 we present moderated versions of the standard likelihood-ratio and Hotelling $T^2$-statistics. We formally build up our multivariate empirical Bayes model and derive the $MB$- and $\tilde{T}^2$-statistics in Section 4. A brief description of a case study is presented in Section 5. Section 6 reports results from a simulation study in which we compare the one-sample $MB$-statistic (the $\tilde{T}^2$-statistic) of Section 4.3 with other statistics. We discuss our models and give directions for future work in Section 7.

Now we introduce some notation for one-sample problems. For two-sample problems, the notation is similar and easily conceived. For each gene $g$, $g = 1, \ldots, G$, we have $n_g$ independent time series, and we model these as i.i.d. $k \times 1$ random vectors $\mathbf{X}_{g1}, \ldots, \mathbf{X}_{gn_g}$ with gene-specific means $\boldsymbol{\mu}_g$ and covariance matrices $\boldsymbol{\Sigma}_g$. Since relative or absolute gene expression measurements are approximately normal on the log scale, we make the multivariate normality assumption on data $\mathbf{X}_{g1}, \ldots, \mathbf{X}_{gn_g}$. Our results are to be judged on their practical usefulness, not on the precise fit of our data to a multivariate normal distribution. As will be seen shortly, our final formulae involve the multivariate $t$ distribution. Thus, a measure of robustness is built in, and our approach will probably be about as effective for elliptically distributed random vectors. We use the natural conjugate priors for $\boldsymbol{\mu}_g$ and $\boldsymbol{\Sigma}_g$, that is, an inverse Wishart prior for $\boldsymbol{\Sigma}_g$ and a dependent multivariate normal prior for $\boldsymbol{\mu}_g$. To simplify the notation, the subscript $g$ will be dropped for the rest of this paper. The statistical models presented in the rest of this paper are for an arbitrary single gene $g$.

The details in this paper differ in two ways from the standard conjugate priors. First, we also have an indicator $I$ such that $I = 1$ when the alternative $K$ is true and $I = 0$ when the null $H$ is true, with the priors for $\boldsymbol{\mu}$ differing in these two cases. Second, when the null hypothesis states that a gene's mean expression level is constant, in order to get a simple closed form expression for the posterior odds, we assume that the gene-specific covariance matrix $\boldsymbol{\Sigma}$ commutes with the $k \times k$ projection matrix $\mathbf{P} = k^{-1}\mathbf{1}_k\mathbf{1}_k'$, that is, $\mathbf{P}\boldsymbol{\Sigma} = \boldsymbol{\Sigma}\mathbf{P}$. In this case the $k \times k$ inverse Wishart prior for $\boldsymbol{\Sigma}$ is replaced by a $(k-1) \times (k-1)$ inverse Wishart prior for a part of $\boldsymbol{\Sigma}$ and an inverse gamma prior for the remainder. These two-part priors are independent; see Section 4.3 for details.

**2. Hypothesis testing.** Our gene ranking problem will be formally stated as a hypothesis testing problem. In this paper we only seek a statistic for ranking genes in the order of evidence against the null hypothesis; we do not hope to obtain raw or adjusted $p$-values as in [27].

Following the notation in [4], the null hypothesis is denoted by $H$, while the alternative hypothesis is denoted by $K$. The null hypothesis corresponding to a gene's mean expression level being 0 is $H: \boldsymbol{\mu} = \mathbf{0}$, $\boldsymbol{\Sigma} > 0$; the alternative is $K: \boldsymbol{\mu} \neq \mathbf{0}$, $\boldsymbol{\Sigma} > 0$. An easy extension is $H: \boldsymbol{\mu} = \boldsymbol{\mu}_0$, $\boldsymbol{\Sigma} > 0$ versus $K: \boldsymbol{\mu} \neq \boldsymbol{\mu}_0$, $\boldsymbol{\Sigma} > 0$, where $\boldsymbol{\mu}_0$ is known. Later we will consider the null hypothesis that a gene's expression is constant, against the alternative that it is not: $H: \boldsymbol{\mu} = \mu_0 \mathbf{1}$, $\boldsymbol{\Sigma} > 0$, where $\mu_0$ is a scalar representing the expected value of the gene's expression level at any time point under $H$, and $\mathbf{1}$ is the $k \times 1$ constant vector of $1s$; $K: \boldsymbol{\mu} \neq \mu_0 \mathbf{1}$, $\boldsymbol{\Sigma} > 0$. Finally, we will consider the null hypothesis that a gene's mean expression levels are the same under two different biological conditions, versus the alternative that they are not: $H: \boldsymbol{\mu}_Z = \boldsymbol{\mu}_Y$, $\boldsymbol{\Sigma}_Z = \boldsymbol{\Sigma}_Y = \boldsymbol{\Sigma} > 0$; $K: \boldsymbol{\mu}_Z \neq \boldsymbol{\mu}_Y$, $\boldsymbol{\Sigma}_Z = \boldsymbol{\Sigma}_Y = \boldsymbol{\Sigma} > 0$.

## 3. The moderated $LR$-statistic.

3.1. *One-sample or paired two-sample problem.* A likelihood-ratio statistic can be used directly to test the null hypothesis $H$ against the alternative hypothesis $K$ when $n > k$. According to standard multivariate results (e.g., [21]), under the alternative hypothesis that there are no constraints on $\boldsymbol{\mu}$ and $\boldsymbol{\Sigma}$, their maximum likelihood estimates are

$$\hat{\boldsymbol{\mu}}_K = \overline{\mathbf{X}},$$

$$\widehat{\boldsymbol{\Sigma}}_K = \frac{n-1}{n} \mathbf{S},$$

where $\mathbf{S} = (n-1)^{-1} \sum_{i=1}^{n} (\mathbf{X}_i - \overline{\mathbf{X}})(\mathbf{X}_i - \overline{\mathbf{X}})'$ is the sample variance–covariance matrix. Also as in [21], the maximum likelihood estimate for the unconstrained $\boldsymbol{\Sigma}$ under the null that $\boldsymbol{\mu} = \boldsymbol{\mu}_0$ is

$$\widehat{\boldsymbol{\Sigma}}_H = \frac{n-1}{n} \mathbf{S} + \mathbf{d}\mathbf{d}',$$

where $\mathbf{d} = \hat{\boldsymbol{\mu}}_K - \hat{\boldsymbol{\mu}}_H$, the difference between the maximum likelihood estimates for $\boldsymbol{\mu}_0$ (if unknown) under $H$ and $K$. If $\boldsymbol{\mu}_0$ is known, then $\mathbf{d} = \hat{\boldsymbol{\mu}}_K - \boldsymbol{\mu}_0$. The likelihood ratio statistic for testing any such $H$ against the above $K$ is

(3.1)
$$\begin{aligned} LR &= 2(l_K^{\max} - l_H^{\max}) \\ &= n \log\left(1 + \frac{n}{n-1} \mathbf{d}' \mathbf{S}^{-1} \mathbf{d}\right). \end{aligned}$$

If the null hypothesis states that $\boldsymbol{\mu}_0 = \mathbf{0}$, then $\mathbf{d}$ reduces to $\overline{\mathbf{X}}$, hence, the likelihood-ratio statistic is equation (3.1) with $\mathbf{d}$ replaced by $\overline{\mathbf{X}}$. Similarly, if $\boldsymbol{\mu}_0$ is known, then $\mathbf{d} = \overline{\mathbf{X}} - \boldsymbol{\mu}_0$. On the other hand, if the null hypothesis states that $\boldsymbol{\mu}_0 = \mu_0 \mathbf{1}$, then the maximum likelihood estimate for $\boldsymbol{\mu}$ is

$$\hat{\boldsymbol{\mu}}_H = \left(\frac{\mathbf{1}' \mathbf{S}^{-1} \overline{\mathbf{X}}}{\mathbf{1}' \mathbf{S}^{-1} \mathbf{1}}\right) \mathbf{1}.$$

The statistic $n\mathbf{d}'\mathbf{S}^{-1}\mathbf{d}$ is the one-sample Hotelling $T^2$-statistic, and by Section 5.3.1b in [21], it is distributed as a Hotelling $T^2(k-1, n-1)$ under $H$, that is, $((n-k+1)n\mathbf{d}'\mathbf{S}^{-1}\mathbf{d})/((n-1)(k-1))$ has an $F$-distribution with degrees of freedom $(k-1, n-k+1)$.

In the microarray time course context, the number of replicates $n$ is typically smaller than the number of time points $k$, and so $\mathbf{S}$ has less than full rank. Furthermore, as discussed in [29], we wish to moderate the sample variance–covariance matrix. Our moderated $\mathbf{S}$ will take the form

$$\widetilde{\mathbf{S}} = \frac{\nu\mathbf{\Lambda} + (n-1)\mathbf{S}}{\nu + n - 1},$$

where $\nu > 0$ controls the degree of moderation, and $\mathbf{\Lambda} > 0$ is the common $k \times k$ matrix toward which $\mathbf{S}$ is smoothed. In Section 4.1 we give the theoretical reason for choosing this moderated variance–covariance matrix $\widetilde{\mathbf{S}}$ and explain how we estimate $\nu$ and $\mathbf{\Lambda}$. Replacing $\mathbf{S}$ with $\widetilde{\mathbf{S}}$ in the $LR$-statistic, our moderated $LR$-statistic is

$$(3.2) \qquad \widetilde{LR} = 2(\boldsymbol{l}_K^{\max} - \boldsymbol{l}_H^{\max}) = n\log\left(1 + \frac{n}{n-1}\tilde{\mathbf{d}}'\widetilde{\mathbf{S}}^{-1}\tilde{\mathbf{d}}\right).$$

When all the genes have an equal number of replicates $n$, equation (3.2) is a monotonic increasing function of $n\tilde{\mathbf{d}}'\widetilde{\mathbf{S}}^{-1}\tilde{\mathbf{d}}$. We define the quadratic form $n\tilde{\mathbf{d}}'\widetilde{\mathbf{S}}^{-1}\tilde{\mathbf{d}} = \|n^{1/2}\widetilde{\mathbf{S}}^{-1/2}\tilde{\mathbf{d}}\|^2$ to be the moderated one-sample Hotelling $T^2$-statistic. In the case of the null $H : \boldsymbol{\mu} = \mathbf{0}, \boldsymbol{\Sigma} > 0$, this is identical to the $\widetilde{T}^2$-statistic we derive in Section 4.1. The one-sample moderated $LR$-statistic and the moderated Hotelling $T^2$-statistic are hybrids of likelihood and Bayesian statistics, since $\widetilde{\mathbf{S}}$ is estimated using the multivariate empirical Bayes procedure we describe below.

3.2. *Unpaired two-sample problem.* Similarly, in the unpaired two-sample case, the moderated $LR$-statistic can be written as a function of the moderated two-sample Hotelling $T^2$-statistic

$$\widetilde{LR} = (m + n)\log\left(1 + \frac{1}{m+n-2}\frac{mn}{m+n}\mathbf{d}'\widetilde{\mathbf{S}}^{-1}\mathbf{d}\right),$$

where $\mathbf{d} = \overline{\mathbf{Z}} - \overline{\mathbf{Y}}$ is the difference between sample averages and $\mathbf{S} = (m + n - 2)^{-1}((m-1)\mathbf{S}_Z + (n-1)\mathbf{S}_Y)$ is the pooled sample variance–covariance matrix, and

$$\widetilde{\mathbf{S}} = \frac{(m + n - 2)\mathbf{S} + \nu\mathbf{\Lambda}}{m + n - 2 + \nu}.$$

The term $(m+n)^{-1}mn\mathbf{d}'\widetilde{\mathbf{S}}^{-1}\mathbf{d}$ is our moderated two-sample Hotelling $T^2$-statistic. We use the same approach to estimate $\widetilde{\mathbf{S}}$ here as that for our two-sample multivariate empirical Bayes model described in Section 4.2.

## 4. The multivariate empirical Bayes model.

### 4.1. *One-sample or paired two-sample problem.*

4.1.1. *Models and priors.* The data $\mathbf{X}_1, \ldots, \mathbf{X}_n$ are multivariate normal with mean $\boldsymbol{\mu}$ and covariance matrix $\boldsymbol{\Sigma}$, denoted by $N_k(\boldsymbol{\mu}, \boldsymbol{\Sigma})$. We define an indicator random variable $I$ to reflect the status of the gene,

$$I = \begin{cases} 1, & \text{if } K \text{ is true,} \\ 0, & \text{if } H \text{ is true.} \end{cases}$$

We suppose that $I$ has a Bernoulli distribution with success probability $p$, $0 < p < 1$. Now we build up our multivariate hierarchical Bayesian model by first assigning independent and identical inverse Wishart priors to the gene-specific covariance matrices $\boldsymbol{\Sigma}$:

$$(4.1) \qquad \boldsymbol{\Sigma} \sim \text{Inv-Wishart}_\nu((\nu\boldsymbol{\Lambda})^{-1}),$$

where $\nu > 0$ and $\nu\boldsymbol{\Lambda} > 0$ are the degrees of freedom and scale matrix, respectively. Given $\boldsymbol{\Sigma}$, we assign multivariate normal priors for the gene-specific mean $\boldsymbol{\mu}$ for the two cases $(I = 1)$ and $(I = 0)$:

$$\boldsymbol{\mu}|\boldsymbol{\Sigma}, I = 1 \sim N_k(\mathbf{0}, \eta^{-1}\boldsymbol{\Sigma}),$$

$$\boldsymbol{\mu}|\boldsymbol{\Sigma}, I = 0 \equiv 0,$$

where $\eta > 0$ is a scale parameter.

The posterior odds are the probability that the expected time course $\boldsymbol{\mu}$ does not stay at 0 (i.e., $I = 1$) over the probability that $\boldsymbol{\mu}$ stays at 0 (i.e., $I = 0$), given the data $\mathbf{X}_1, \ldots, \mathbf{X}_n$. Following [24]'s notation, we write

$$(4.2) \qquad \begin{aligned} \mathbf{O} &= \frac{P(I = 1|\mathbf{X}_1, \ldots, \mathbf{X}_n)}{P(I = 0|\mathbf{X}_1, \ldots, \mathbf{X}_n)} \\ &= \frac{p}{1 - p} \frac{P(\mathbf{X}_1, \ldots, \mathbf{X}_n|I = 1)}{P(\mathbf{X}_1, \ldots, \mathbf{X}_n|I = 0)}. \end{aligned}$$

The distribution of the data given $I$ can be written as

$$(4.3) \qquad P(\mathbf{X}_1, \ldots, \mathbf{X}_n|I) = \int P(\overline{\mathbf{X}}|\boldsymbol{\Sigma}, I) P(\mathbf{S}|\boldsymbol{\Sigma}, I) P(\boldsymbol{\Sigma}|I) \, d\boldsymbol{\Sigma}.$$

4.1.2. *Multivariate joint distributions.* Once the priors and the models are set, the joint distributions of the data can be determined given $I$. We omit the standard calculations leading to

$$\begin{aligned} & P(\mathbf{X}_1, \ldots, \mathbf{X}_n|I = 1) \\ &= \frac{\boldsymbol{\Gamma}_k((n + \nu)/2)}{\boldsymbol{\Gamma}_k((n - 1)/2)\boldsymbol{\Gamma}_k(\nu/2)} \end{aligned}$$

(4.4)

$$\times (n-1)^{(1/2)k(n-1)} \nu^{-(1/2)kn} \pi^{-(1/2)k} (n^{-1} + \eta^{-1})^{-(1/2)k}$$

$$\times \frac{|\boldsymbol{\Lambda}|^{-(1/2)n} |\mathbf{S}|^{(1/2)(n-k-2)}}{|\mathbf{I}_k + ((n^{-1} + \eta^{-1})\nu\boldsymbol{\Lambda})^{-1}\overline{\mathbf{X}\mathbf{X}}' + (\nu\boldsymbol{\Lambda}/(n-1))^{-1}\mathbf{S}|^{(1/2)(n+\nu)}}.$$

Thus, given $I = 1$, the probability density function of the data is a function of $\overline{\mathbf{X}}$ and $\mathbf{S}$ only, which follows a Student–Siegel distribution [1]. Following [1]'s notation, this distribution is denoted by $StSi_k(\nu, \mathbf{0}, (n^{-1} + \eta^{-1})\boldsymbol{\Lambda}, n - 1, (n-1)^{-1}\nu\boldsymbol{\Lambda})$. For the case of $I = 0$, we get the same distribution with different parameters, namely, $StSi_k(\nu, \mathbf{0}, n^{-1}\boldsymbol{\Lambda}, n - 1, (n-1)^{-1}\nu\boldsymbol{\Lambda})$.

4.1.3. *MB-statistic and $\widetilde{T}^2$-statistic.* Define our moderated gene-specific sample variance–covariance matrix $\widetilde{\mathbf{S}}$ to be the inverse of the posterior mean of $\boldsymbol{\Sigma}^{-1}$ given $\mathbf{S}$,

(4.5)
$$\widetilde{\mathbf{S}} = [E(\boldsymbol{\Sigma}^{-1}|\mathbf{S})]^{-1} = \frac{(n-1)\mathbf{S} + \nu\boldsymbol{\Lambda}}{n-1+\nu}.$$

The posterior odds $\mathbf{O}$ we defined earlier can be derived using the distributions of the data given $I$ and is

(4.6)
$$\mathbf{O} = \frac{p}{1-p} \left(\frac{\eta}{n+\eta}\right)^{(1/2)k}$$
$$\times \left(\frac{n-1+\nu+\widetilde{T}^2}{n-1+\nu+(\eta/(n+\eta))\widetilde{T}^2}\right)^{(1/2)(n+\nu)},$$

where $\widetilde{T}^2 = \widetilde{\mathbf{t}}'\widetilde{\mathbf{t}}$ and $\widetilde{\mathbf{t}}$ is the moderated multivariate $t$-statistic defined by

(4.7)
$$\widetilde{\mathbf{t}} = n^{1/2}\widetilde{\mathbf{S}}^{-1/2}\overline{\mathbf{X}}.$$

Following the tradition in genetics, the log base 10 of $\mathbf{O}$ is called the LOD score. To distinguish it from the LOD score (also called the $B$-statistic) in the univariate model of [20] and [24], the multivariate LOD score in this paper is called the *MB-statistic*,

(4.8)
$$MB = \log_{10} \mathbf{O}.$$

When all genes have the same number of replicates $n$, equation (4.8) is a monotonic increasing function of $\widetilde{T}^2$. This shows that the *MB*-statistic is equivalent to the $\widetilde{T}^2$-statistic when $n$ is the same across genes, and therefore, one is encouraged to use the $\widetilde{T}^2$-statistic in this case since it does not require the estimation of $\eta$ and leads to the same rankings as equation (4.8). We now derive the distribution for $\widetilde{T}^2$.

By Gupta and Nagar [16], the Jacobian transformation from $\overline{\mathbf{X}}$ to $\widetilde{\mathbf{t}}$ is $J(\overline{\mathbf{X}} \to \widetilde{\mathbf{t}}) = |n^{-1/2}\widetilde{\mathbf{S}}^{1/2}|$. Since equation (4.4) is a function of $\overline{\mathbf{X}}$ and $\mathbf{S}$ only,

it is the joint probability density function for these two random variables. Substituting for $\overline{\mathbf{X}}$ in terms of $\tilde{\mathbf{t}}$ in equation (4.4), and multiplying the resulting expression by $J(\overline{\mathbf{X}} \to \tilde{\mathbf{t}})$, the joint probability density function for $\tilde{\mathbf{t}}$ and $\mathbf{S}$ given $I = 1$ is

$$
\begin{aligned}
P(&\tilde{\mathbf{t}}, \mathbf{S} | I = 1) \\
&= \pi^{-(1/2)k} \frac{\Gamma((n+\nu)/2)}{\Gamma((n+\nu-k)/2)} \left( \frac{n+\eta}{\eta} \right)^{-(1/2)k} (n-1+\nu)^{-(1/2)k} \\
&\quad \times \left( 1 + \frac{1}{n-1+\nu} \left( \frac{\eta}{n+\eta} \right) \tilde{\mathbf{t}}'\tilde{\mathbf{t}} \right)^{-(1/2)(n+\nu)} \\
&\quad \times \frac{1}{\boldsymbol{\beta}_k((n-1)/2, \nu/2)} \\
&\quad \times \frac{|\mathbf{S}|^{(1/2)(n-k-2)}}{|\nu\boldsymbol{\Lambda}/(n-1)|^{(1/2)(n-1)} |\mathbf{I}_k + (\nu\boldsymbol{\Lambda}/(n-1))^{-1}\mathbf{S}|^{(1/2)(n+\nu-1)}}.
\end{aligned}
$$
(4.9)

The above expression is factorized into parts involving $\mathbf{S}$ only and $\tilde{\mathbf{t}}$ only, proving that $\tilde{\mathbf{t}}$ and $\mathbf{S}$ are independent. It is apparent that $\tilde{\mathbf{t}}$ has a multivariate $t$ distribution with degrees of freedom $n+\nu-k$, scale parameter $n+\nu-1$, mean vector $\mathbf{0}$ and covariance matrix $\eta^{-1}(n+\eta)\mathbf{I}_k$. This distribution is denoted by $\tilde{\mathbf{t}} | I = 1 \sim \mathbf{t}_k(n+\nu-k, n+\nu-1, \mathbf{0}, \eta^{-1}(n+\eta)\mathbf{I}_k)$ [16]. It is straightforward to see that $\tilde{\mathbf{t}} | I = 0 \sim \mathbf{t}_k(n+\nu-k, n+\nu-1, \mathbf{0}, \mathbf{I}_k)$. Given $I = 1$, $\mathbf{S}$ is distributed as a generalized type-II beta distribution with parameters $(n-1)/2$, $\nu/2$, scale matrix $\nu\boldsymbol{\Lambda}/(n-1)$ and location matrix $\mathbf{0}$. The distribution is denoted by $GB_k^{II}((n-1)/2, \nu/2, \nu\boldsymbol{\Lambda}/(n-1), \mathbf{0})$ [16]. The marginal distribution of $\mathbf{S}$ does not depend on $I$ so that $P(\mathbf{S}|I=0) = P(\mathbf{S}|I=1)$. This distributional result is used to estimate the hyperparameter $\boldsymbol{\Lambda}$. The distribution for $\tilde{T}^2$ under the null follows immediately. Under $H$, $k^{-1}\tilde{T}^2$ has an $F$ distribution with degrees of freedom $(k, n+\nu-k)$; equivalently, $(n+\nu-k)^{-1}(n+\nu-1)\tilde{T}^2$ has the Hotelling $T^2$-distribution $T^2(k, n+\nu-1)$.

The $\tilde{T}^2$-statistic is identical to the one-sample moderated Hotelling $T^2$-statistic in Section 3.1 with the same null hypothesis.

For the easy extension to the above model, $H : \boldsymbol{\mu} = \boldsymbol{\mu}_0, \boldsymbol{\Sigma} > 0$ and $K : \boldsymbol{\mu} \neq \boldsymbol{\mu}_0, \boldsymbol{\Sigma} > 0$, where $\boldsymbol{\mu}_0$ is known, all the results above hold with $\overline{\mathbf{X}}$ replaced by $\overline{\mathbf{X}} - \boldsymbol{\mu}_0$.

### 4.1.4. *Special cases.*

1. $\boldsymbol{\Sigma} = \sigma^2 \mathbf{I}_k$.   By constraining $\boldsymbol{\Sigma} = \sigma^2 \mathbf{I}_k$, we ignore the correlations among gene expression values at different times, and assume the variances at different times are equal. Suppose that the prior for $\sigma^2$ is

$$
\sigma^2 \sim \text{inv-gamma}(\tfrac{1}{2}\nu, \tfrac{1}{2}\nu\lambda^2).
$$

Define

$$s_j^2 = (n-1)^{-1} \sum_{i=1}^n (\mathbf{X}_{ij} - \overline{\mathbf{X}}_j)^2,$$

$$\tilde{s}_j^2 = (n-1+\nu)^{-1}((n-1)s_j^2 + \nu\lambda^2)$$

and

$$\tilde{t}_j = n^{1/2}\overline{\mathbf{X}}_j \tilde{s}_j^{-1}, \qquad j = 1,\dots,k.$$

In this case, the posterior odds are equivalent to a product of $k$ independent univariate odds,

$$\mathbf{O} = \frac{p}{1-p}\left(\frac{\eta}{n+\eta}\right)^{(1/2)k} \prod_{j=1}^k \left(\frac{n-1+\nu+\tilde{t}_j^2}{n-1+\nu+(\eta/(n+\eta))\tilde{t}_j^2}\right)^{(1/2)(n+\nu)},$$

(4.10)

and the *MB*-statistic is equivalent to the sum of $k$ univariate *B*-statistics.

2. $n = 1$. When $n = 1$, that is, when there is no replication at all, each gene has its own unknown variability. The moderated multivariate $t$-statistic becomes $\tilde{\mathbf{t}} = \mathbf{\Lambda}^{-1/2}\mathbf{X}$. The posterior odds are obtained by plugging in $n = 1$ in equation (4.6), and are found to be a function of $\mathbf{X}$ only. Since there is no replication, our hyperparameters must be assigned values, for example, from previous experiments.

3. $k = 1$. When $k = 1$, that is, when there is only one time point, the alternative hypothesis states that there is differential expression at this single time point. Our multivariate model should and does reduce to the univariate model in [20] and [24].

### 4.1.5. *Limiting cases.*

1. $\nu \to \infty$. In this case, the gene-specific variance–covariance matrices are totally ignored. The moderated multivariate $t$-statistic above becomes $\tilde{\mathbf{t}}_\infty = n^{1/2}\mathbf{\Lambda}^{-1/2}\overline{\mathbf{X}}$, and $\widetilde{T}_\infty^2 = \tilde{\mathbf{t}}_\infty'\tilde{\mathbf{t}}_\infty$. The posterior odds become

$$\mathbf{O} = \frac{p}{1-p}\left(\frac{\eta}{n+\eta}\right)^{(1/2)k} \exp\left(\frac{1}{2}\left(\frac{n}{n+\eta}\right)\widetilde{T}_\infty^2\right).$$

2. $\nu \to 0$. In this case, there is no moderation at all. The posterior odds are just equation (4.6) with $\nu$ replaced by 0. If $n < k$, then $\mathbf{S}^{-1/2}$ should be calculated by a $g$-inverse.

3. $\nu \to \infty$ *and* $\boldsymbol{\Sigma} = \sigma^2 \mathbf{I}_k$. Define $\tilde{t}_{\infty j} = n^{1/2} \lambda^{-1} \overline{\mathbf{X}}_j$, $j = 1, \ldots, k$. The posterior odds become

$$\mathbf{O} = \frac{p}{1-p} \left( \frac{\eta}{n+\eta} \right)^{(1/2)k} \exp\left( \frac{1}{2} \left( \frac{n}{n+\eta} \right) \sum_{j=1}^{k} \tilde{t}_{\infty j}^2 \right).$$

4. $\nu \to 0$ *and* $\boldsymbol{\Sigma} = \sigma^2 \mathbf{I}_k$. In this case, the posterior odds are just equation (4.10) with $\nu$ replaced with 0.

4.1.6. *Hyperparameter estimation.* We have shown that the *MB*-statistic for assessing whether or not a time course has mean $\mathbf{0}$ depends on $(k^2 + k + 6)/2$ hyperparameters: $\nu$, $\boldsymbol{\Lambda}$, $\eta$ and $p$. In practice, we need to estimate these hyperparameters, and plug in our estimates into the formulae for $\bar{\mathbf{S}}$, $\tilde{\mathbf{t}}$, $\mathbf{O}, \ldots$ and so on. Slightly abusing our notation, we will use the same symbols for these estimates, relying on context to make it clear whether we are assuming the hyperparameters to be known or not. In our multivariate model, many more hyperparameters need to be estimated, compared to the univariate models in [20] and [24], both of which have four hyperparameters. Closed form estimators for the hyperparameters in the univariate linear model setting are derived in [24], using the marginal sampling distributions of the statistic $\tilde{t}$ and the sample variance $s^2$, and are shown to be better than the simple estimators in [20]. Following [24], the aim of this section is to derive estimators for the hyperparameters in our multivariate model. In general, the hyperparameter $\eta$ associated with the case $I = 1$ is estimated based on only a small subset of genes, while $\nu$ and $\boldsymbol{\Lambda}$ are estimated using the whole gene set. Instead of estimating the proportion of differentially expressed genes $p$, we plug in a user-defined value, since the choice of $p$ does not affect the rankings of genes based on the *MB*-statistic.

*EB estimation of $\nu$ and $\boldsymbol{\Lambda}$.* The hyperparameter $\nu$ determines the degree of smoothing between $\mathbf{S}$ and $\boldsymbol{\Lambda}$. The method we use to estimate $\nu$ builds on that used to estimate $d_0$ in Section 6.2 in [24]. However, unlike $d_0$ in [24], $\nu$ is associated with the $k \times k$ matrix $\boldsymbol{\Sigma}$. Therefore, a method appropriate to this multivariate framework is needed. Let $\hat{\nu}_j$ be the estimated prior degrees of freedom based on the $j$th diagonal elements of the gene-specific sample variance–covariance matrices (i.e., the replicate variances for the $j$th time point from the whole gene set) using the method proposed in Section 6.2 in [24]. Our estimation of $\nu$ is based on the following two-step strategy. As the first step, set $\nu$ as $\hat{\nu} = \max(\text{mean}(\hat{\nu}_j), k + 6)$, $j = 1, \ldots, k$. This estimated $\hat{\nu}$ is used to estimate $\boldsymbol{\Lambda}$. Once $\boldsymbol{\Lambda}$ is estimated, $\hat{\nu}$ is reset to be $\hat{\nu} = \text{mean}(\hat{\nu}_j)$. In practice, one can even just plug in a user-defined value $\nu_0$ which gives the desired amount of smoothing. In such a case, the first step sets $\hat{\nu} =$ $\max(\nu_0, k + 6)$. This $\hat{\nu}$ is used to estimate $\boldsymbol{\Lambda}$. After $\boldsymbol{\Lambda}$ is estimated, $\hat{\nu}$ can be reset to the user-defined value $\nu_0$.

Our estimate of $\boldsymbol{\Lambda}$ comes after the first step in the estimation of $\nu$. We showed that, under our model, $\mathbf{S}$ follows the generalized type-II beta distribution with expectation $(\nu - k - 1)^{-1}\nu\boldsymbol{\Lambda}$. By the weak law of large numbers, $\overline{\mathbf{S}}$ converges in probability to $(\nu - k - 1)^{-1}\nu\boldsymbol{\Lambda}$. We can thus estimate $\boldsymbol{\Lambda}$ by $\hat{\nu}^{-1}(\hat{\nu} - k - 1)\overline{\mathbf{S}}$. If $\hat{\nu} \to \infty$, then $\boldsymbol{\Lambda}$ is estimated by $\overline{\mathbf{S}}$. The above estimates give quite satisfactory results on real data. A theoretical analysis of the estimation of our hyperparameters will be given later. For the moment, we content ourselves with obtaining reasonable estimates.

*EB estimation of $\eta$.* The hyperparameter $\eta$ is related to the moderated multivariate $t$-statistic $\tilde{\mathbf{t}}$ of nonzero genes. The method we use to estimate $\eta$ builds on that of estimating $v_0$ in [24], except that we now need to deal with the multivariate case. Let $\tilde{t}_j$ be the $j$th element of the moderated multivariate $t$-statistic $\tilde{\mathbf{t}}$, $j = 1, \ldots, k$. As in Section 6.3 in [24], each $\tilde{t}_j$ gives an estimate of $\eta$, call it $\hat{\eta}_j$, based on the top $p/2$ portion of genes with the largest $|\tilde{t}_j|$. We set $\hat{\eta}$ to be the mean of the $\hat{\eta}_j$.

4.2. *Unpaired two-sample problem.* Suppose there are two independent biological conditions $Z$ and $Y$ with sample sizes $m$ and $n$, respectively. We can also derive the *MB*-statistic for testing the null $H : \boldsymbol{\mu}_Z = \boldsymbol{\mu}_Y, \boldsymbol{\Sigma}_Z = \boldsymbol{\Sigma}_Y > 0$. The null hypothesis turns out to be the same as that in Section 4.1: $H : \boldsymbol{\mu} = \mathbf{0}, \boldsymbol{\Sigma} > 0$, if we write $\boldsymbol{\mu} = \boldsymbol{\mu}_Z - \boldsymbol{\mu}_Y$ and $\boldsymbol{\Sigma} = \boldsymbol{\Sigma}_Z = \boldsymbol{\Sigma}_Y$. That is, we solve this two-sample problem using the one-sample approach in Section 4.1. We denote the $m$ i.i.d. time course vectors for biological condition $Z$ by $\mathbf{Z}_1, \ldots, \mathbf{Z}_m$, each from a multivariate normal distribution with mean $\boldsymbol{\mu}_Z$ and variance–covariance matrix $\boldsymbol{\Sigma}$. Similarly, those for biological condition $Y$ are denoted by $\mathbf{Y}_1, \ldots, \mathbf{Y}_n$, each with mean $\boldsymbol{\mu}_Y$ and variance–covariance matrix $\boldsymbol{\Sigma}$. Since the null hypothesis here is identical to that in Section 4.1, the priors for $\boldsymbol{\mu}$ and $\boldsymbol{\Sigma}$ are exactly the same as those in Section 4.1, and we omit the details here. In a later paper we will attack this problem by assigning independent priors for $\boldsymbol{\mu}_Y$ and $\boldsymbol{\mu}_Z$ separately.

All the results follow immediately. The moderated multivariate $t$-statistic $\tilde{\mathbf{t}}$ here is defined as equation (4.7) with $n$ replaced by $(m^{-1} + n^{-1})^{-1}$ and $\overline{\mathbf{X}}$ replaced by $\overline{\mathbf{Z}} - \overline{\mathbf{Y}}$. $\tilde{\mathbf{S}}$ here is the same as that defined in Section 3.2. The posterior odds $\mathbf{O}$ against the null hypothesis that the expected time courses are the same are

$$\mathbf{O} = \frac{p}{1-p}\left(\frac{m^{-1} + n^{-1}}{m^{-1} + n^{-1} + \eta^{-1}}\right)^{(1/2)k}$$
$$\times \left(\frac{m + n - 2 + \nu + \tilde{T}^2}{m + n - 2 + \nu + ((m^{-1} + n^{-1})/(m^{-1} + n^{-1} + \eta^{-1}))\tilde{T}^2}\right)^{(1/2)(m+n+\nu-1)}.$$

The log base 10 of $\mathbf{O}$ is our two-sample $MB$-statistic. Again, when all genes have the same sample sizes $m$ and $n$, the two-sample $MB$-statistic is equivalent to the $\tilde{T}^2 = \tilde{\mathbf{t}}' \tilde{\mathbf{t}}$. Under $H$, $k^{-1} \tilde{T}^2$ has an $F$ distribution with degrees of freedom $(k, m+n+\nu-k-1)$; equivalently, $(m+n+\nu-k-1)^{-1}(m+n+\nu-2) \tilde{T}^2$ has the Hotelling $T^2$ distribution $T^2(k, m+n+\nu-2)$.

The $MB$-statistic described in this section involves hyperparameters $\nu$, $\mathbf{\Lambda}$, $\eta$ and $p$. The estimation procedures for these hyperparameters are very similar to those in Section 4.1, except that we have to use the gene-specific pooled sample variance–covariance matrices when estimating $\nu$ and $\mathbf{\Lambda}$, and use the $k \times 1$ moderated multivariate $t$-statistic $\tilde{\mathbf{t}}$ defined here to estimate $\eta$. We omit the details here.

It should be noted that the $MB$-statistic derived in this section has a slightly different definition: instead of using all the data we observe, we only use the difference in sample averages and the pooled sample variance–covariance matrix. The $\tilde{T}^2$-statistic here is identical to the moderated two-sample Hotelling $T^2$-statistic in Section 3.2.

4.3. *One-sample problem of constancy.* In this section we derive the posterior odds against the null that a gene's mean expression level stays constant over time. We obtain a closed form solution similar to that in the preceding sections, but only under a constraint on the variance–covariance matrix $\mathbf{\Sigma}$.

4.3.1. *Transformation.* For each gene, let $I$ be the indicator variable defined in Section 4.1. Let $\mathbf{P} = k^{-1} \mathbf{1}_k \mathbf{1}'_k$ be the $k \times k$ projection matrix onto the rank 1 space of constant vectors, where $\mathbf{1}'_k = (1, \ldots, 1)$ is a $k \times 1$ vector of 1s. Let $\mathbf{P}^c = \mathbf{I}_k - \mathbf{P}$ be the projection onto the orthogonal complement of $R(\mathbf{P})$. We can write any vector $\boldsymbol{\mu} \in R^k$ as $\boldsymbol{\mu} = \mathbf{P}\boldsymbol{\mu} + \mathbf{P}^c\boldsymbol{\mu}$, and in the case $I = 0$, the second term $\mathbf{P}^c\boldsymbol{\mu}$ vanishes. As in Section 4.1, we build up our multivariate model by first assigning independent inverse Wishart priors to the gene-specific covariance matrices $\mathbf{\Sigma}$; see equation (4.1). Given $\mathbf{\Sigma}$, we next assign multivariate normal priors to the gene-specific mean parameters $\boldsymbol{\mu}$ for the case of nonconstant $(I = 1)$ and constant genes $(I = 0)$, respectively:

$$(4.11) \qquad \begin{cases} \boldsymbol{\mu}|\mathbf{\Sigma}, I = 1 \sim N(0, \tau^{-1}\mathbf{P}\mathbf{\Sigma}\mathbf{P} + \kappa^{-1}\mathbf{P}^c\mathbf{\Sigma}\mathbf{P}^c), \\ \boldsymbol{\mu}|\mathbf{\Sigma}, I = 0 \sim N(0, \tau^{-1}\mathbf{P}\mathbf{\Sigma}\mathbf{P}). \end{cases}$$

Given $\mathbf{\Sigma}$ and $I = 0$, the covariance matrix $\mathbf{P}\mathbf{\Sigma}\mathbf{P}$ guarantees that $\boldsymbol{\mu}$ is a constant vector, while when $I = 1$, the extra component $\mathbf{P}^c\mathbf{\Sigma}\mathbf{P}^c$ adds further variance to $\boldsymbol{\mu}$ so that it becomes a nonconstant vector. Again, in order to obtain the full expression for the posterior odds $\mathbf{O}$, we need to derive $P(\mathbf{X}_1, \ldots, \mathbf{X}_n | I)$ using equation (4.3). To get a closed-form expression for the posterior odds, we find it necessary to make an additional assumption, namely, that $\mathbf{P}\mathbf{\Sigma} = \mathbf{\Sigma}\mathbf{P}$. With this assumption, given $\mathbf{\Sigma}$ and $I = 0$, $\overline{\mathbf{X}}$ is a multivariate normal distribution with mean $\mathbf{0}$ and covariance matrix $(n^{-1}\mathbf{\Sigma} + \tau^{-1}\mathbf{\Sigma P})$. Similarly, given $\mathbf{\Sigma}$ and $I = 1$, $\overline{\mathbf{X}}$ is a multivariate normal distribution with mean $\mathbf{0}$ and the covariance matrix $(n^{-1}\mathbf{\Sigma} + \tau^{-1}\mathbf{\Sigma P} + \kappa^{-1}\mathbf{\Sigma P}^c)$.

For the rest of this section, unless stated otherwise, we assume $\mathbf{P\Sigma} = \mathbf{\Sigma P}$, and we make use of the following lemma, whose proof is omitted.

LEMMA 4.1. *Suppose* $\mathbf{T}$ *is any* $k \times k$ *nonsingular matrix whose first row is constant* $c$ *and the remaining rows have row sums equal to 0. Write* $\mathbf{T} = (\mathbf{T}_0', \mathbf{T}_1')'$, *where* $\mathbf{T}_0$ *is the first row of* $\mathbf{T}$, *and* $\mathbf{T}_1$ *is the remainder. Then, for any* $\mathbf{\Sigma} > 0$ *satisfying* $\mathbf{P\Sigma} = \mathbf{\Sigma P}$, $\mathbf{T\Sigma T}' = \widetilde{\mathbf{\Sigma}}$ *is a* $k \times k$ *block diagonal matrix with the scalar* $\tilde{\sigma}^2 > 0$ *as the first block and* $(k-1) \times (k-1)$ *matrix* $\widetilde{\mathbf{\Sigma}}_1 > 0$ *as the second block: that is,*

$$\mathbf{T\Sigma T}' = \widetilde{\mathbf{\Sigma}} = \begin{pmatrix} \tilde{\sigma}^2 & \mathbf{0} \\ \mathbf{0} & \widetilde{\mathbf{\Sigma}}_1 \end{pmatrix}.$$

As the first example, let $\mathbf{T}$ be the Helmert matrix, where the $ji$th element of $\mathbf{T}$ is defined as

$$\begin{cases} t_{ji} = 1/\sqrt{k}, & \text{for } j = 1, i = 1, \dots, k, \\ t_{ji} = 1/\sqrt{j(j-1)}, & \text{for } 2 \le j \le k, 1 \le i \le j-1, \\ t_{ji} = -(j-1)/\sqrt{j(j-1)}, & \text{for } 2 \le j \le k, i = j, \\ t_{ji} = 0, & \text{for } 2 \le j \le k-1, j+1 \le i \le k. \end{cases}$$

$\mathbf{T}$ can also be the following matrix, where the $ji$th element of $\mathbf{T}$ is defined as

$$\begin{cases} t_{ji} = 1, & \text{for } j = 1, i = 1, \dots, k, \\ t_{ji} = 1, & \text{for } 2 \le j \le k, i = j-1, \\ t_{ji} = -1, & \text{for } 2 \le j \le k, i = j, \\ t_{ji} = 0, & \text{otherwise.} \end{cases}$$

For our multivariate empirical Bayes model in this section, we use the Helmert matrix $\mathbf{T}$ to proceed with our calculations. The results can be applied to other $\mathbf{T}$ immediately.

4.3.2. *Models and priors.* Here $\mathbf{T}$ is partitioned into its first row $\mathbf{T}_0$ $(1 \times k)$ and its last $k-1$ rows $\mathbf{T}_1$ $((k-1) \times k)$. Since $\mathbf{X}_1, \dots, \mathbf{X}_n$ are i.i.d. $N_k(\boldsymbol{\mu}, \mathbf{\Sigma})$, the transformed random vectors $\mathbf{TX}_i$ are also multivariate normally distributed with mean $\mathbf{T}\boldsymbol{\mu}$ and covariance matrix $\widetilde{\mathbf{\Sigma}}$. By Lemma 4.1, the matrix $\widetilde{\mathbf{\Sigma}}$ is a block diagonal matrix with $\tilde{\sigma}^2$ as the first block, and $\widetilde{\mathbf{\Sigma}}_1$ as the second block. Defining $\bar{x}_i = k^{-1} \sum_{j=1}^k X_{ij}$, then $\sqrt{k}\bar{x}_i$ and the random vectors $\mathbf{T}_1\mathbf{X}_i$ are independent and normally distributed, with distributions

$$\sqrt{k}\bar{x}_i | \mathbf{T}_0\mu, \tilde{\sigma}^2 \sim N(\mathbf{T}_0\boldsymbol{\mu}, \tilde{\sigma}^2),$$

$$\mathbf{T}_1\mathbf{X}_i | \mathbf{T}_1\boldsymbol{\mu}, \widetilde{\mathbf{\Sigma}}_1 \sim N(\mathbf{T}_1\boldsymbol{\mu}, \widetilde{\mathbf{\Sigma}}_1).$$

This transformation allows us to separate the gene expression changes into constant and nonconstant changes.

As we have seen in Section 4.1, the joint distributions of data given $I$ can be fully described using the sufficient statistics $\bar{\bar{x}}$, $\mathbf{T}_1\overline{\mathbf{X}}$, $s^2$ and $\mathbf{S}_1$, where $\bar{\bar{x}} = n^{-1}\sum_{i=1}^{n}\bar{x}_i$, $\mathbf{T}_1\overline{\mathbf{X}} = n^{-1}\sum_{i=1}^{n}\mathbf{T}_1\mathbf{X}_i$, $s^2 = (n-1)^{-1}\sum_{i=1}^{n}(\bar{x}_i - \bar{\bar{x}})^2$ and $\mathbf{S}_1 = (n-1)^{-1}\sum_{i=1}^{n}(\mathbf{T}_1\mathbf{X}_i - \mathbf{T}_1\overline{\mathbf{X}})(\mathbf{T}_1\mathbf{X}_i - \mathbf{T}_1\overline{\mathbf{X}})'$. The prior for $\widetilde{\mathbf{\Sigma}}$ is first set through the independent priors for $\tilde{\sigma}^2$ and $\widetilde{\mathbf{\Sigma}}_1$. We suppose that $\tilde{\sigma}^2$ and $\widetilde{\mathbf{\Sigma}}_1$ are independently distributed, with an inverse gamma distribution with shape parameter $\xi/2$ and scale parameter $\xi\lambda^2/2$, and an inverse Wishart distribution with degrees of freedom $\nu$ and scale matrix $\nu\mathbf{\Lambda}_1$, respectively, that is,

$$(4.12) \qquad \begin{cases} \tilde{\sigma}^2 \sim \text{inv-gamma}(\frac{1}{2}\xi, \frac{1}{2}\xi\lambda^2), \\ \widetilde{\mathbf{\Sigma}}_1 \sim \text{inv-Wishart}_\nu((\nu\mathbf{\Lambda}_1)^{-1}). \end{cases}$$

The prior for $\mathbf{T}\boldsymbol{\mu}$ has four parts. We assign independent priors to $\mathbf{T}_0\boldsymbol{\mu}$ and $\mathbf{T}_1\boldsymbol{\mu}$ separately for the cases $I = 1$ and $I = 0$. For the case $I = 1$, priors are

$$(4.13) \qquad \begin{cases} \mathbf{T}_0\boldsymbol{\mu}|\tilde{\sigma}^2, I = 1 \sim N(\theta, \kappa^{-1}\tilde{\sigma}^2), \\ \mathbf{T}_1\boldsymbol{\mu}|\widetilde{\mathbf{\Sigma}}_1, I = 1 \sim N(\mathbf{0}, \eta^{-1}\widetilde{\mathbf{\Sigma}}_1), \end{cases}$$

where $\theta \geq 0$ is the mean, and $\kappa > 0$ and $\eta > 0$ are scale parameters. When $I = 0$, $\mathbf{T}_1\boldsymbol{\mu} = \mathbf{0}$ with probability 1. Thus, the priors in this case are

$$(4.14) \qquad \begin{cases} \mathbf{T}_0\boldsymbol{\mu}|\tilde{\sigma}^2, I = 0 \sim N(\theta, \kappa^{-1}\tilde{\sigma}^2), \\ \mathbf{T}_1\boldsymbol{\mu}|\widetilde{\mathbf{\Sigma}}_1, I = 0 \equiv \mathbf{0}. \end{cases}$$

It is reasonable to assume $P(\mathbf{T}_0\boldsymbol{\mu}|\tilde{\sigma}^2, I = 0) = P(\mathbf{T}_0\boldsymbol{\mu}|\tilde{\sigma}^2, I = 1)$ for large genome-wide arrays since there is no obvious reason why the expected grand mean of the expression levels for nonconstant genes should differ from that of constant genes. For two-color comparative microarray experiments, it is also reasonable to assume $\theta = 0$.

4.3.3. *Multivariate joint distributions.* The joint distributions can be derived quite readily using a previously established formula, and so we omit the details. Given $I = 1$, $\mathbf{T}_1\overline{\mathbf{X}}$ and $\mathbf{S}_1$ follow the Student–Siegel distribution $StSi_{k-1}(\nu, \mathbf{0}, (n^{-1} + \eta^{-1})\mathbf{\Lambda}_1, n-1, (n-1)^{-1}\nu\mathbf{\Lambda}_1)$. Similarly, the joint distribution of $\mathbf{T}_1\overline{\mathbf{X}}$ and $\mathbf{S}_1$ given $I = 0$ is $StSi_{k-1}(\nu, \mathbf{0}, n^{-1}\mathbf{\Lambda}_1, n-1, (n-1)^{-1}\nu\mathbf{\Lambda}_1)$.

4.3.4. *MB-statistic and $\widetilde{T}^2$-statistic.* The posterior odds against the null that a gene's mean expression level stays constant over time are equation (4.6) in Section 4.1 with $k$ replaced by $k-1$, $\tilde{\mathbf{t}}$ expressed by equation (4.7) with $\widetilde{\mathbf{S}}$ replaced by $\widetilde{\mathbf{S}}_1$ and $\overline{\mathbf{X}}$ replaced by $\mathbf{T}_1\overline{\mathbf{X}}$. $\widetilde{\mathbf{S}}_1$ is just equation (4.5) with $\mathbf{S}$ replaced by $\mathbf{S}_1$ and $\mathbf{\Lambda}$ replaced by $\mathbf{\Lambda}_1$. As in Section 4.1, the *MB*-statistic is a monotonic increasing function of $\widetilde{T}^2 = \widetilde{\mathbf{t}}'\widetilde{\mathbf{t}}$ when all genes have the same sample size $n$.

Under $H$, $(k-1)^{-1}\widetilde{T}^2$ has an $F$ distribution with degrees of freedom $(k-1, n+\nu-k+1)$, or, equivalently, $(n+\nu-k+1)^{-1}(n+\nu-1)\widetilde{T}^2$ has a Hotelling $T^2$-distribution $T^2(k-1, n+\nu-1)$. The hyperparameter estimation procedures here are similar to those described in Section 4.1, except that all the estimations are performed based on transformed data.

**5. Case study.** In this section we illustrate our results with a paired two-sample problem, using the *Arabidopsis thaliana* dataset in [35]. Here we only give a very brief description of the data and the results. We refer the reader to [35] and to Chapter 5 of [28] for more thorough discussions.

*A. thaliana* wildtype (Columbia) and ics1-2 null mutant plants were evenly positioned, intermixed and grown in growth chambers under controlled conditions. When the plants were four weeks old, they were infected with a moderately heavy innoculum of the powdery mildew *G. orontii*. Each pair of mRNA samples from wildtype and mutant plants was harvested and collected at six time points post-infection. Plants could not be resampled, so mRNA samples at one time point were from different plants than those of any other time point. We report here on the analysis of three replicate experiments under similar environmental conditions which contribute four biological replicates: one from the first and third experiments (1–3, 3–1) and two from the second experiment (2–1, 2–2). These mRNA samples were hybridized onto Affymetrix *Arabidopsis* ATH1 GeneChips, yielding 22,810 probesets and 48 arrays in our analysis. The array preprocessing were done using the Robust Multi-array Analysis (RMA) algorithm described in [5, 18] which is implemented in the Bioconductor package `affy`.

This study is longitudinal if we treat *experiments* as units, while it is cross-sectional if we treat *plants* as units. We believe it is worthwhile to treat it as a paired longitudinal study, since samples within the same experiment are more similar than those from different experiments. We thus have a paired two-sample problem, with the genes of interest being those whose wildtype and mutant temporal profiles are different. We subtracted the $\log_2$ intensities of the wildtype from those of the paired mutant at each time within each replicate, yielding the $\log_2$ ratios of mutant relative to wildtype. Since the number of time course replicates is the same ($n = 4$) across genes for this dataset, we used the $\widetilde{T}^2$-statistic instead of the *MB*-statistic to rank genes, so that we did not have to estimate the hyperparameter $\eta$.

For comparison purposes, we fitted a linear model to the log ratios for each gene, with time and replicate effects, and calculated the $F$-statistic for the time effect.

TABLE 1
*Spearman rank correlations between $\widetilde{T}^2$ with different $\nu$ and the estimated $\nu$ for all and the top 859 genes. The percent moderation is defined by $(\nu + n - 1)^{-1}\nu \times 100$*

| % moderation | Correlation (all) | Correlation (top 859) |
|---|---|---|
| 97 ($\nu = 100$) | 0.97 | 0.90 |
| 80 ($\nu = 12$) | 0.99 | 0.98 |
| 40 ($\nu = 2$) | 0.99 | 0.98 |
| 25 ($\nu = 1$) | 0.98 | 0.96 |
| 0 ($\nu = 0.01$) | 0.93 | 0.90 |

5.1. *Results.* The extent of moderation from $\hat{\nu} = 5$ was 63%. The left panel of Figure 1 presents three genes falling into different ranges of ranks (rank = 1, 175, 859) by $\widetilde{T}^2$, while the ones on the right panel have the same ranks by $F$. The gene ranked most highly by $\widetilde{T}^2$ exhibits much greater differences between the wildtype and mutant temporal profiles than the one ranked most highly by $F$. The magnitude of the difference, as measured by $\widetilde{T}^2$, decreases as the rank goes down. The gene ranked 1 by $\widetilde{T}^2$ is well known: pathogenesis-related protein 1 (PR1). Other known pathogenesis-related genes also ranked highly by $\widetilde{T}^2$ and were less highly ranked by the $F$-statistic, as detailed in [35].

To investigate the effect of the amount of moderation on gene ranking, we kept $\mathbf{\Lambda}$ fixed, and re-calculated the $\widetilde{T}^2$-statistic with several different $\nu$ values. We then computed the Spearman rank correlation between the different sets of $\widetilde{T}^2$-statistics for all genes and for the top 859 genes only. Table 1 gives the results. The correlations are lower in the two extremes. All of these sets have the same number one gene. This comparison shows that the gene ranks are reasonably stable when the extent of moderation varies within a certain window, and that moderation seems to have more effect on the top genes relative to the whole gene set.

## 6. Simulation study.

6.1. *Method.* In this section we report on a small simulation study for the null hypothesis $H : \boldsymbol{\mu} = \mu_0 \mathbf{1}$, $\mathbf{\Sigma} > 0$ based on an actual example we have met. We simulate 100 data sets, each with 20,000 genes. The genes are simulated independently, which we regard as an assumption that makes sense to compare methods, but it should be kept in mind that gene expression measures in real data can be quite dependent. In each simulated data set, 400 out of the 20,000 genes are assigned to be nonconstant. That is, $p = 0.02$. Each gene is simulated with three independent replicates ($n = 3$) and eight time points ($k = 8$). The other hyperparameters are the following: $\nu = 13$, $\xi = 3$, $\lambda^2 = 0.3$, $\theta = 0$ (two-color experiments), $\kappa = 0.02$, $\eta = 0.08$, and

$$\boldsymbol{\Lambda} = \begin{pmatrix} 14.69 & 0.57 & 0.99 & 0.40 & 0.55 & 0.51 & -0.23 \\ 0.57 & 15.36 & 1.22 & 0.84 & 1.19 & 0.91 & 0.86 \\ 0.99 & 1.22 & 14.41 & 2.47 & 1.81 & 1.51 & 1.07 \\ 0.40 & 0.84 & 2.47 & 17.05 & 2.40 & 2.32 & 1.33 \\ 0.55 & 1.19 & 1.81 & 2.40 & 15.63 & 3.31 & 2.75 \\ 0.51 & 0.91 & 1.51 & 2.32 & 3.31 & 13.38 & 3.15 \\ -0.23 & 0.86 & 1.07 & 1.33 & 2.75 & 3.15 & 12.90 \end{pmatrix} \times 10^{-3}.$$

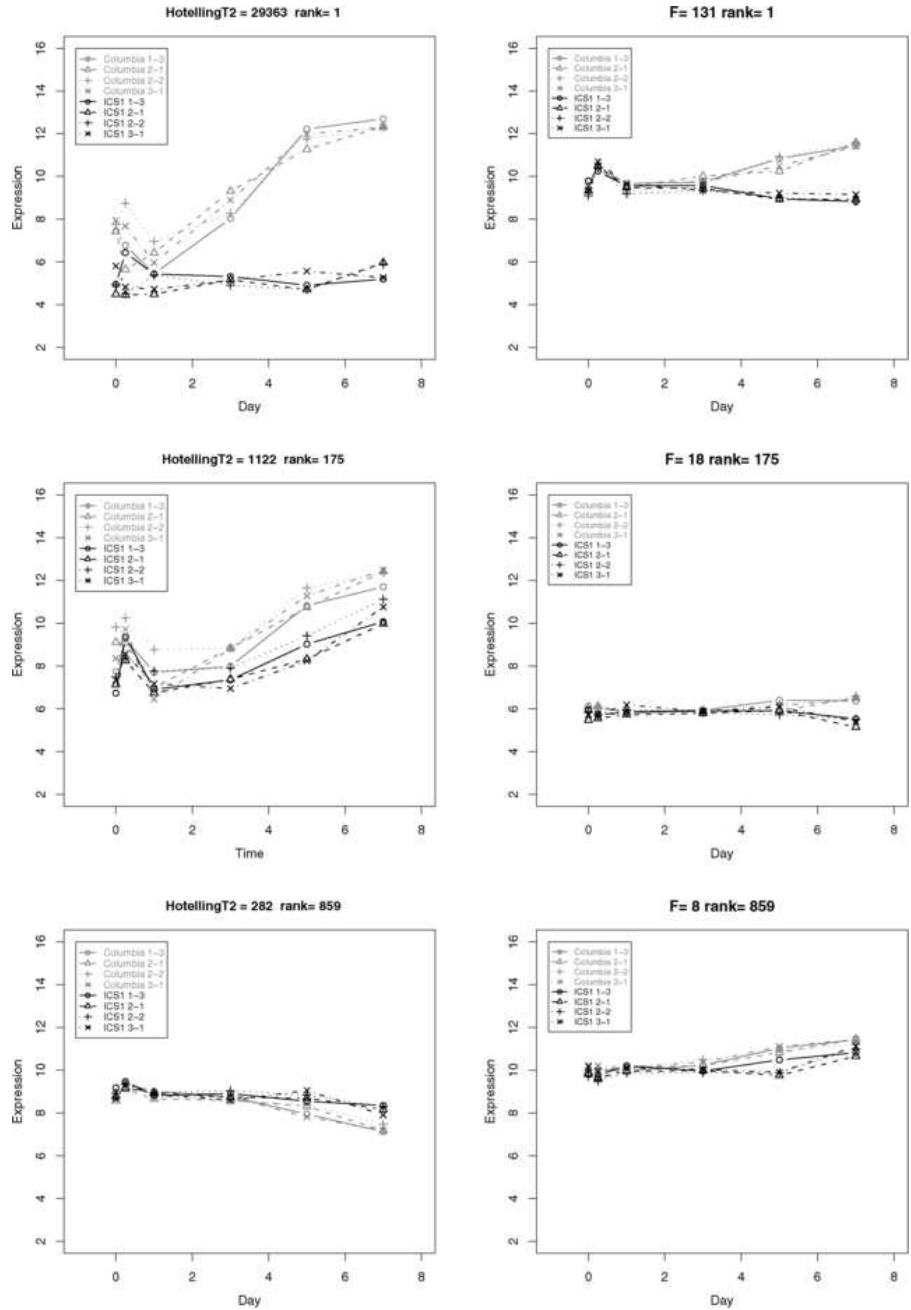

The correlation matrix of $\boldsymbol{\Lambda}$ is

$$\begin{pmatrix} 1 & 0.04 & 0.07 & 0.03 & 0.04 & 0.04 & -0.02 \\ 0.04 & 1 & 0.08 & 0.05 & 0.08 & 0.06 & 0.06 \\ 0.07 & 0.08 & 1 & 0.16 & 0.12 & 0.11 & 0.08 \\ 0.03 & 0.05 & 0.16 & 1 & 0.15 & 0.15 & 0.09 \\ 0.04 & 0.08 & 0.12 & 0.15 & 1 & 0.23 & 0.20 \\ 0.04 & 0.06 & 0.11 & 0.15 & 0.23 & 1 & 0.24 \\ -0.02 & 0.06 & 0.08 & 0.09 & 0.19 & 0.24 & 1 \end{pmatrix},$$

and we see clear evidence of serial correlation. Note that, in the real world, although often the case, the correlation does not always decrease with time lag. The statistics compared in our study are the following:

(1) *MB*-statistic, or equivalently, the $\widetilde{T}^2$-statistic;

(2) *MB*-statistic using first differences: take the differences in gene expression values at consecutive time points within replicates, and use them to test the null hypothesis $H : \boldsymbol{\mu} = \mathbf{0}, \boldsymbol{\Sigma} > 0$ (Section 4.1), where $\boldsymbol{\mu}$ is the mean of the differences;

(3) *MB*-statistic in the special case $\boldsymbol{\Sigma} = \sigma^2 \mathbf{I}_k$;

(4) *MB*-statistic in the limiting case $\nu \to \infty$;

(5) *MB*-statistic in the limiting case $\nu \to 0$;

(6) ordinary $F$-statistic from an ANOVA model with time and replicate effects;

(7) partly-moderated $F$-statistic proposed in [24] from an ANOVA model with time and replicate effects;

(8) one-sample moderated Hotelling $T^2$-statistic $\|n^{1/2} \widetilde{\mathbf{S}}^{-1/2} \tilde{\mathbf{d}}\|^2$ derived in Section 3, equivalently, the moderated *LR*-statistic, where the degree of moderation and the common matrix toward which each sample covariance matrix moves is estimated by the method given in Section 4.1;

(9) the variance across time course replicates $(nk-1)^{-1} \sum_{i=1}^{n} \sum_{j=1}^{k} (X_{ij} - \bar{\bar{x}})^2$.

Here each of the nine statistics incorporates either none (e.g., variance) or one (ordinary $F$-statistic) or more of the following: moderation, correlation structure and replicate variances, and thus can be used to show the importance of the above properties. It is not appropriate to set the prior degrees of freedom $\nu$ to be a very small number, since we have the constraint that $\nu \geq k - 1$. We choose $\nu$ to be $k + 5 = 13$ because it simulates more stable $\mathbf{\Sigma}$'s across genes.

TABLE 2
*The means and standard deviations (SD) of the diagonal elements of the estimated $\mathbf{\Lambda}_1$*

| Hyperparameters | True value $\times 10^3$ | Mean $\times 10^3$ | SD $\times 10^3$ |
|---|---|---|---|
| $\lambda_1^2$ | 14.69 | 14.71 | 0.16 |
| $\lambda_2^2$ | 15.36 | 15.37 | 0.17 |
| $\lambda_3^2$ | 14.41 | 14.43 | 0.15 |
| $\lambda_4^2$ | 17.05 | 17.04 | 0.19 |
| $\lambda_5^2$ | 15.63 | 15.63 | 0.15 |
| $\lambda_6^2$ | 13.38 | 13.40 | 0.15 |
| $\lambda_7^2$ | 12.90 | 12.92 | 0.17 |

6.2. *Results.* Table 2 compares the means and standard deviations of the hyperparameter estimates of the diagonal elements of $\mathbf{\Lambda}_1$ ($\lambda_j^2$), $j = 1, \ldots, k-1$, with their true values. The mean estimate of $\mathbf{\Lambda}_1$ is very close to the true $\mathbf{\Lambda}_1$, and the standard deviations are very small. The hyperparameter $\eta$ is always under-estimated (mean $= 0.026$, SD $= 0.002$), which agrees with Section 8 in [24], where $v_0$ was usually over-estimated. The hyperparameter $\nu$ is also always under-estimated (mean $= 7.0$, SD $= 0.2$). In Section 5 we observed that the amount of moderation $\nu$ does not greatly affect gene ranking except at the two extremes. One can even choose a user-defined $\nu$ which gives reasonable results. Although not well estimated, $\eta$ only affects the rankings when the number of replicates is different across genes. However, this does not happen often in the real world. Even when that happens, the effect is very small. To investigate the effects of $\eta$ on gene rankings, we tried a couple of $\eta$ values from different ranges, while keeping the remaining hyperparameters fixed, and calculated the *MB*-statistics. The rank correlations between rankings of the *MB*-statistics with the user-defined $\eta$'s and the estimated $\eta$ for one set of simulated data are the following: 0.91, 0.94, 0.99, 0.99, 0.99 for $\eta = 2$, 1, 0.08 (true value), 0.05, 0.001, respectively.

To examine the relationship between the $\widetilde{T}^2$-statistic and the true deviation from constancy, the $\log_{10}$ transformed $\widetilde{T}^2$-statistic from one simulated dataset is plotted against the Mahalanobis distance between the expected time course vector $\boldsymbol{\mu}$ and its projection onto the rank 1 constant space $\bar{\boldsymbol{\mu}} = \mathbf{P}\boldsymbol{\mu}$ (Figure 2). The squared Mahalanobis distance is defined by $d(\boldsymbol{\mu}, \bar{\boldsymbol{\mu}})^2 = (\boldsymbol{\mu} - \bar{\boldsymbol{\mu}})' \boldsymbol{\Sigma}^{-1} (\boldsymbol{\mu} - \bar{\boldsymbol{\mu}})$. Figure 2 clearly shows that $\log_{10} \widetilde{T}^2$ is positively correlated with $d(\boldsymbol{\mu}, \bar{\boldsymbol{\mu}})$, and most of the 400 true nonconstant genes achieve higher $\widetilde{T}^2$-statistics than the constant genes.

Figure 3 plots the average number of false positives against average number of false negatives at different cutoffs. By different cutoffs, we mean choosing the top $x$ genes and calculating the numbers of false positives and false negatives, where $x$ varies across the integers from 400 to 800. The lines in Figure 3 from left to right represent the following: *MB*-statistic ($\widetilde{T}^2$), *MB*-statistic with first differences (indistinguishable from the *MB*-statistic), one-sample moderated Hotelling $T^2$-statistic (indistinguishable from the *MB*-statistic), *MB*-statistic with $\boldsymbol{\Sigma} = \sigma^2 \mathbf{I}_k$, *MB*-statistic with $\nu \to \infty$, partly-moderated $F$-statistic [24], ordinary $F$-statistic, *MB*-statistic with $\nu \to 0$ and variance. The *MB*-statistic ($\widetilde{T}^2$) attains almost the same numbers of false positives and false negatives as *MB* with first differences and the one-sample moderated Hotelling $T^2$-statistic. The effectiveness of moderation is demonstrated by comparing the lines for the *MB*-statistic, the *MB*-statistic in the limiting case $\nu \to \infty$ and the *MB*-statistic in the limiting case that $\nu \to 0$. Both of these limiting cases lead to higher aggregate false positives and false negatives. This result supports the view stated in [29] that moderation is useful. In particular, the case $\nu \to 0$ (no moderation at all) produces much higher numbers of false positives and false negatives. This is likely due to the poor estimation of sample variance–covariance matrices with a small number of replicates. Indeed, the ordinary unmoderated $F$-statistic which ignores the correlation structure achieves smaller numbers of false positives and false negatives than the unmoderated $MB$-statistic. A similar situation also arises in the microarray discrimination context; see Section 7 of [12]. The partly-moderated $F$-statistic [24] which ignores the dependence among times behaves like the $MB$-statistic in the special case $\mathbf{\Sigma} = \sigma^2 \mathbf{I}_k$. Moreover, it achieves fewer false positives and false negatives than the ordinary $F$-statistic. Figure 3 also demonstrates the importance of incorporating the correlation structure among time points. The $MB$-statistic ($\tilde{T}^2$) and the one-sample moderated Hotelling $T^2$-statistic perform better than the partly-moderated $F$-statistic in [24] and the ordinary $F$-statistic; the former incorporates the correlation structure among time points, whereas the latter does not. However, we observe that the amount of moderation given by the partly-moderated $F$-statistic in [24] is usually much less than that given by the $MB$-statistic. When there are a large number of residual degrees of freedom from the linear model, the partly-moderated $F$-statistic [24] behaves very much like the ordinary $F$-statistic. This suggests that the lower number of false positives and number of false negatives from the $MB$-statistic than the partly-moderated $F$-statistic in [24] involve both the incorporation of correlation structures and the amounts of moderation. As expected, the simple variance statistic across replicates, which totally ig-nores the replicate variances, performs the worst. This demonstrates the importance of incorporating the replicate variances into any statistic.

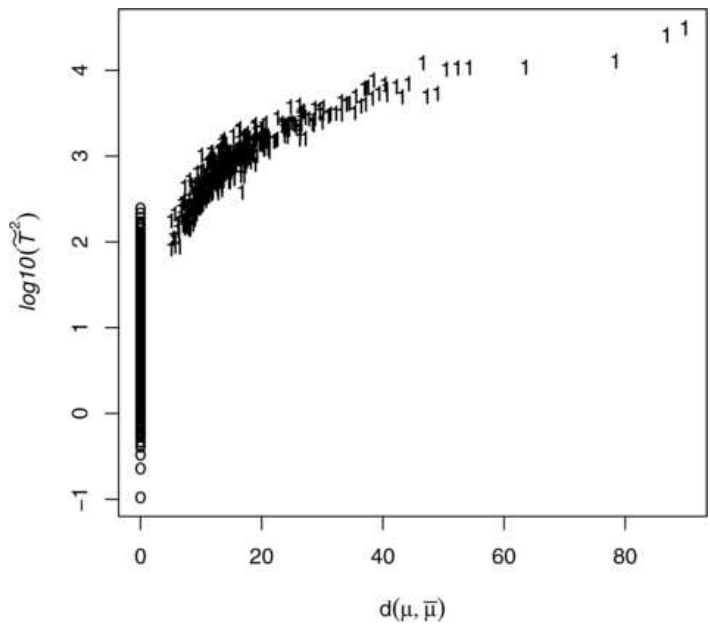

Fig. 2. *The* $\log_{10} \widetilde{T}^2$ *statistic versus the true deviation from constancy* $d(\boldsymbol{\mu}, \bar{\boldsymbol{\mu}})$ *for one simulated dataset. Here* 1 *denotes nonconstant, and* o *constant genes.*

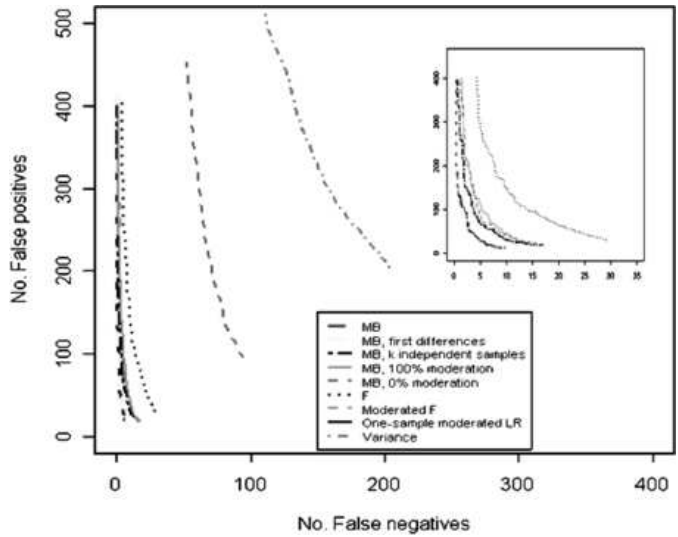

Fig. 3. *Average number of false positives versus number of false negatives of all the nine statistics. The subplot presents the curves for the best seven statistics.*

**7. Discussion.** In this paper we introduced the $MB$- and $\widetilde{T}^2$-statistics for one- and two-sample longitudinal replicated developmental microarray time course experiments. Our main focus was the one-sample or paired two-sample problem with the null hypothesis $H : \boldsymbol{\mu} = \mathbf{0}, \boldsymbol{\Sigma} > \mathbf{0}$. This $MB$-statistic can be used when there are two biological conditions and the samples are paired across conditions, and it is shown to perform better than the classical $F$-statistic on a problem briefly described in Section 5. In addition, we also derive the $MB$-statistics and $\widetilde{T}^2$-statistics for the two-sample problem with the null $H : \boldsymbol{\mu}_Z = \boldsymbol{\mu}_Y, \boldsymbol{\Sigma}_Z = \boldsymbol{\Sigma}_Y = \boldsymbol{\Sigma} > 0$, and the one-sample problem with the null $H : \boldsymbol{\mu} = \mu_0 \mathbf{1}, \boldsymbol{\Sigma} > \mathbf{0}$ using similar approaches. The latter situation requires a slight assumption on $\boldsymbol{\Sigma}$ in order to get a simple closed-form solution for the posterior odds against the null. All the $MB$-statistics and $\widetilde{T}^2$-statistics incorporate the correlation structure, replication and moderation. The moderated versions of some standard likelihood-ratio test statistics are also described. When all genes have the same sample size(s), our $\widetilde{T}^2$-statistics are not only equivalent to the $MB$-statistics, but also are identical to their corresponding moderated Hotelling $T^2$-statistics, apart from the one-sample problem with the null $H : \boldsymbol{\mu} = \mu_0 \mathbf{1}, \boldsymbol{\Sigma} > \mathbf{0}$, where there is an additional constraint on $\boldsymbol{\Sigma}$. In this case the $\widetilde{T}^2$-statistic performed as well as the moderated Hotelling $T^2$-statistic in our simulation study, and also on several real datasets we have encountered. We have shown in the simulation study that, with this null, the $MB$-statistic ($\widetilde{T}^2$), the $MB$-statistic using first differences and the one-sample moderated Hotelling $T^2$-statistic perform best among all the nine statistics compared. This is not entirely surprising given that we simulated data under our model, but the comparisons are still informative. In practice, we regard the $MB$-statistic, the $MB$-statistic with first differences and the moderated $LR$-statistic as performing equally well, giving very similar (or identical) results. However, to use the $LR$-statistic (the moderated Hotelling $T^2$-statistic), we still need to insert moderated sample variance–covariance matrices, and these come from our multivariate empirical Bayes framework. In other words, our models provide a natural way to estimate the gene-specific moderated variance–covariance matrices (Sections 4.1–4.3), while the likelihood-based approach alone does not.

The assumption of $\mathbf{P}\boldsymbol{\Sigma} = \boldsymbol{\Sigma}\mathbf{P}$ with the null $H : \boldsymbol{\mu} = \mu_0 \mathbf{1}, \boldsymbol{\Sigma} > 0$ allows the mathematical calculations in Section 4.3, and leads to a closed-form formula for the $MB$-statistic. One question which naturally arises to be the impact of this constraint on the rankings of genes. From the practical point of view, the impact of this constraint on gene rankings is very *slight*. The rank correlations between the one-sample $MB$-statistic *with* the commuting assumption and the moderated $LR$-statistic from likelihood-based approach *without* the constraint, from the actual examples we have met, are typically very high (over 0.99). The rank correlations from our simulated data are also over 0.99. On the other hand, using the $MB$-statistic with first differences produces very similar or even identical results to the $MB$-statistic in Section 4.3. Indeed, instead of using the Helmert matrix, if we choose $\mathbf{T}$ to be the second transformation matrix of Section 4.3, we get identical results. Even so, we do not consider the first differences approach to be the solution to this commuting constraint, since the inference drawn is based on reduced, not the original data. In other words, the null hypotheses are not equivalent. The likelihood-based approach with moderation and the first differences approach support the fact that this constraint does not have much effect on the results. In addition, the former is a good way to avoid this assumption, since it performs as well as our $MB$-statistic in Section 4.3.

The statistics proposed in this paper are for one- and two-sample longitudinal data. We should be aware that many experiments in the real world exhibit some features from both longitudinal and cross-sectional experiments (e.g., Section 5).

One thing we plan to investigate in the future is the effect of assuming the same variance–covariance matrix $\mathbf{\Sigma}$ for both $I = 1$ and $I = 0$. Another issue which interests us is the effect of assuming the same $\mathbf{\Sigma}$ across biological conditions in the unpaired two-sample model in Section 4.2. The proposed methods may be extended in several ways, for example, identifying genes of some *specific* pattern, rather than *any* pattern. The statistics for a longitudinal multi-sample problem when there are at least three biological conditions and genes of interest are those with different temporal profiles across conditions derived in [30]. The corresponding statistics for cross-sectional data are also presented in [31].

The proposed methods focus on gene ranking, but not assessing the *significance* using $p$-values. However, if desired, we believe that generating $p$-values from a bootstrap analysis should be successful in this context.

We constructed our models using conjugate priors for the multivariate normal likelihoods, so that we got simple closed-form solutions for the posteriors odds. Finding a closed-form statistic when the priors on $\boldsymbol{\mu}$ and $\mathbf{\Sigma}$ are independent seems to be an open and probably hard problem; that problem probably needs to be dealt with using MCMC.

**Acknowledgments.** We thank the Editor Jianqing Fan, an Associate Editor and the two referees for their valuable comments on this paper. We are grateful to Mary Wildermuth and her colleagues for sharing their *A. thaliana* dataset prior to publication of Wildermuth et al. [35]. Thanks are also due to Gordon Smyth, Ingrid Lönnstedt, Darlene Goldstein, Greg Hather, Avner

Bar-Hen, Alfred Hero, Cavan Reilly and Christina Kendziorski for their valuable comments on an earlier draft, and to Ben Bolstad for his assistance with our simulations. Finally, we would like to acknowledge Mary Wildermuth, John Ngai and their lab members for their helpful discussions on the biological background of time course experiments and access to their data to test the methods developed in this paper.

Division of Biostatistics
367 Evans Hall 3860
University of California, Berkeley
Berkeley, California 94720-3860
USA
E-mail: yuchuan@stat.berkeley.edu
        terry@stat.berkeley.edu